\documentclass[12pt]{article}
\usepackage{amsmath}
\usepackage{amssymb,color}

\def\0{{\bf 0}}
\def\1{{\bf 1}}
\def\l{{\ell}}

\def\proof{\noindent{\bf Proof: }}
\def\qed{ \hskip 20pt{\vrule height7pt width6pt depth0pt}\hfil}
\def\forb{{\mathrm{forb}}}
\def\ext{{\mathrm{ext}}}
\def\Av{{\mathrm{Avoid}}}

\newcommand{\linelessfrac}[2]{\genfrac{}{}{0pt}{}{#1}{#2}}

\newcommand{\ncols}[1]{\| #1 \|}
\newcommand{\rf}[1]{(\ref{#1})}
\newcommand{\trf}[1]{Theorem~\ref{#1}}

\newcommand{\lrf}[1]{Lemma~\ref{#1}}

\newtheorem{thm}{Theorem}[section]
\newtheorem{lemma}[thm]{Lemma}

\newtheorem{conj}[thm]{Conjecture}

\newenvironment{E}{\begin{equation}}{\end{equation}}
\usepackage[english]{babel}
\makeatother

\usepackage{graphicx} 
\title{
   Forbidden Configurations and Boundary Cases}
\author{R.P. Anstee\thanks{Research supported in part by
NSERC}, Oakley Edens\thanks{Research supported in part by NSERC USRA}, Arvin Sahami\thanks{Research supported in part by NSERC USRA}, \\ Mathematics Department\\The University of British Columbia\\Vancouver,
B.C. Canada V6T 1Z2\\ {\small{\texttt{anstee@math.ubc.ca}}}, 
{\small{\texttt{oedens@math.harvard.edu}}},{\small{\texttt{sahamiarvin52@gmail.com}}},\\
\and Attila Sali\thanks{Research  partially supported by the
    National Research, Development and Innovation Office (NKFIH)
    grants K--132696 and SNN-135643. 
}\\ HUN-REN Alfr\'ed R\' enyi Institute of Mathematics\\Budapest, Hungary and \\ Department of Computer Science\\ Budapest University of Technology and Economics\\ {\small{\texttt{sali.attila@renyi.hu}}} 
} 

\begin{document}

\maketitle

\begin{abstract}
    Let $F$ be a $k\times \l$ (0,1)-matrix. Define a (0,1)-matrix $A$ to have a $F$ as a \emph{configuration} if there is a submatrix of $A$ which is a row and column permutation of $F$. In the language of sets, a configuration is a \emph{trace}.
Define a matrix to be {\it simple} if it is a (0,1)-matrix with no repeated columns. Let $\Av(m,F)$ be all simple $m$-rowed matrices $A$ with no configuration $F$.  Define $\forb(m,F)$ as the maximum number of columns of any matrix in $\Av(m,F)$.   Determining $\forb(m,F)$ requires determining bounds and constructions of matrices in $\Av(m,F)$. 
The paper  considers some column maximal $k$-rowed simple $F$ that have the bound $\Theta(m^{k-2})$ and yet adding a column increases bound to $\Omega(m^{k-1})$. By a construction,  $\forb(m,F)$ is determined exactly. 
 \hfil\break
Keywords: extremal set theory,  (0,1)-matrices, forbidden configurations, trace, constant weight codes
\end{abstract}

\section{Introduction}

The paper considers some boundary forbidden configurations. We start with some notation.

An $m\times n$ matrix $A$ is said to be \emph{simple} if it is a (0,1)-matrix with no repeated columns. There is a natural correspondence between columns of $A$ and subsets of $[m]$. We consider an extremal set problem in matrix terminology as follows. Let $\ncols{A}$ be the number of columns of $A$. For a given matrix $F$, we say $F$ is a \emph{configuration} in $A$, denoted $F\prec A$, if there is a submatrix of $A$ which is a row and column permutation of $F$. Define 
$$\Av(m,F)=\left\{A\,|\,A \hbox{ is }m\hbox{-rowed and simple}, F\not\prec A\right\},$$ 
$$\forb(m,F)=\max_{A\in\Av(m,F)}\ncols{A}.$$
A matrix $A\in\Av(m,F)$ is called {\emph{extremal}} if $\ncols{A}=\forb(m,F)$ and let 
$$\ext(m,F)=\{A\in\Av(m,F)\,|\,\ncols{A}=\forb(m,F)\}.$$ 

Many results are recorded in a survey\cite{survey}.  \trf{mainthm} obtains exact bounds for four extremal configurations, cases of $k$-rowed simple $F$ where $\forb(m,F)$ is $O(m^{k-2})$ and adding any column $\alpha$ to $F$ has $\forb(m,[F \,|\,\alpha])$ being $\Omega(m^{k-1})$. The main  conjecture in \cite{survey} suggests which matrices have this property. We discuss this after  
\trf{simple,k-2}.

 Define $f(m,k)$ for $k\ge 2$ by the recurrence 
 \begin{E}f(m,k)=f(m-1,k)+f(m-1,k-1),\label{recurrence}\end{E}
\noindent with the base cases  $f(m,2)=2$ and $f(2,k)=2$ for $k>2$ (so that \break$f(m,3)=2m$ as in \trf{F1}).
  
  We can solve the recurrence to obtain by induction
$$f(m,k)=f(m-1,k)+f(m-1,k-1)\hfil$$
$$=\binom{m-2}{k-2}+\sum_{i=0}^{k-2}\binom{m-1}{i}+\binom{m-2}{k-3}+\sum_{i=0}^{k-3}\binom{m-1}{i}
=\binom{m-1}{k-2}+ \sum_{i=0}^{k-2}\binom{m}{i}.$$
Note $f(m,2)=2$. Also 
\begin{E}f(m,k)=2\sum_{i=0}^{k-2}\binom{m-1}{i}
 =2\binom{m}{k-2}+2\binom{m}{k-4}+2\binom{m}{k-6}+\cdots \label{f(m,k)} \end{E}.

\begin{thm}\cite{survey} 
\cite{AFl10}\label{simple,k-2} Let $k\ge 2$ be given.

\noindent If ${\cal F}$ is a family of simple $k\times \l$ matrices with the property that there is an $F\in{\cal F}$ with a pair of rows that do not contain $K_2^0$, that there is an $F\in{\cal F}$ with a pair of rows that do not contain $K_2^2$ and that there is an $F\in{\cal F}$ with a pair of rows that do not contain the configuration $K_2^1=I_2$, then $\forb(m,F)$ is $O(m^{k-2})$, more precisely
$\forb(m,{\cal F})\le f(m,k)$ defined by the recurrence \rf{recurrence} with bases cases $f(m,2)=2$ and $f(2,k)=2$ for $k>2$ so that $f(m,k)=\binom{m-1}{k-2}+\sum_{i=0}^{k-2}\binom{m}{i}$.
\\ \noindent If $F$ is a simple $k\times \l$ matrix with the property that either every pair of rows has $K_2^0$ or every pair of rows has $K_2^2$ or every pair of rows has $K_2^1$, then $\forb(m,F)$ is $\Theta(m^{k-1})$.\qed 
\end{thm}
\vskip 10pt
Some important matrices include $I_k$, the $k\times k$ identity matrix, $T_k$ the $k\times k$ triangular matrix (with $1's$ in position $i,j$ if $i\le j$), and $K_k$ the $k\times 2^k$ matrix of all possible (0,1)-columns on $k$ rows. We denote by $F^c$ the (0,1)-complement of $F$ so that $I_k^c$ is the complement of the identity. Define $K_k^s$ to be the $k\times \binom{k}{s}$ matrix of columns of sum  $s$. We define $\1_k$ as $k\times 1$ column of 1's, $\0_k$ as $k\times 1$ column of 0's and $\1_k\0_{\ell}$ as $(k+\ell)\times 1$ columns with $k$ 1's on top of $\ell$ 0's.  

For two simple matrices $A,B$ where $A$ is $m_1$-rowed and $B$ is $m_2$-rowed, the \emph{product} $A\times B$ is defined as the $(m_1+m_2)$-rowed matrix of $\ncols{A}\ncols{B}$ columns consisting of each column of $A$ on top of each column of $B$. For example $K_k=[0\,1]\times K_{k-1}$. 
For a subset $S$ of rows, define $A_S$ as the submatrix of $A$ formed by those rows.  
We typically ignore row and column permutations of our matrices unless explicitly stated.

The conjecture \cite{survey} asserts that the asymptotic bounds are achieved by products of $I,I^c,T$. 

\begin{conj}Let $F$ be given. Let $t$ be the largest integer so that there is a $t$-fold product of matrices chosen from $I,I^c,T$ that avoids $F$. Then $\forb(m,F)$ is $\Theta(m^t)$.\end{conj}
Applying this allows one to
identify possible \emph{boundary cases}, namely configurations $F$, with conjectured bound $\forb(m,F)$, where  adding any column to $F$ forming $F'$ has $\forb(m,F')$ being $\Omega(m\cdot\forb(m,F))$.    All simple $k$-rowed $F$, for which $\forb(m,F)$ is $O(m^{k-2})$ yet adding any column increases bound to $\Omega(m^{k-1})$, is given below in \rf{6boundary}.  The $(k$-$2)$-fold product construction $A=I_{m/(k-2)}\times I_{m/(k-2)}\times\cdots\times  I_{m/(k-2)}$ has $A\in\Av(m,F)$. In particular the column maximal property ensures that there is only one pair of rows of $F$ avoiding $\1_2$.

$$F_{1,k}=\begin{array}{c}
\left[\begin{array}{cccc}1&1&1&0\\ 1&0&0&1\\ 0&1&0&0\\ \end{array}\right]\\
\times \\
K_{k-3}\\ \end{array}
,\quad
F_{2,k}=\begin{array}{c}
 \left[\begin{array}{cc}1&0\\ 0&1\\ \end{array}\right]\\
\times \\
\left[\begin{array}{ccc}0&1&1\\ 0&0&1\\ \end{array}\right]\\
\times\\
K_{k-4}\\ \end{array},$$

\begin{E}F_{3,k}=\begin{array}{c}
\left[\begin{array}{cccccccc} 
0&0&0&0&1&0&0&1\\ 
0&0&0&1&0&1&1&0\\
0&1&1&0&0&1&1&1\\
1&0&1&1&1&0&1&1\\
\end{array}\right]\\
\times\\ 
K_{k-4}\\ \end{array},
\quad
F_{4,k}=\begin{array}{c}
\left[\begin{array}{ccccc} 
0&0&1&1&0\\        
0&0&0&0&1\\
0&1&0&1&1\\
\end{array}\right]\\
\times\\
\left[\begin{array}{ccc}
1&0&1\\  
0&1&1\\  
 \end{array}\right]\\
\times\\
K_{k-5}\\ \end{array},\label{6boundary}\end{E}

$$F_{5,k}=\begin{array}{@{}c@{}}
\left[\begin{array}{@{}ccccccccccccccc@{}} 
1&1&1&0&0&0&0&0&0&1&1&1&1&1&1\\        
0&0&0&1&1&1&1&1&1&1&1&1&1&1&1\\
0&0&0&0&0&0&1&1&1&0&0&0&1&1&1\\
0&1&0&0&1&0&0&1&0&0&1&0&0&1&0\\
0&0&1&0&0&1&0&0&1&0&0&1&0&0&1\\
\end{array}\right]\\
\times\\
K_{k-5}\\ \end{array} 
,\,\,
F_{6,k}= \begin{array}{@{}c@{}}
\left[\begin{array}{@{}ccc@{}}0&1&0\\ 0&0&1\\ \end{array}\right]\\
\times \\
\left[\begin{array}{@{}ccc@{}}1&0&1\\ 0&1&1\\ \end{array}\right]\\
\times \\ 
\left[\begin{array}{@{}ccc@{}}0&1&1\\ 0&0&1\\ \end{array}\right]\\
\times\\
K_{k-6}\\ \end{array}$$
This paper establishes some exact bounds as follows.

\begin{thm} $\forb(m,F_{1,4})=f(m,4)$ and $\forb(m,F_{3,k})=f(m,k)$ for $k=4,5,6$. \label{mainthm} \end{thm}

Note $f(m,k)$ is $\Theta(m^{k-2})$.  Now $A=I_{m/(k-2)}\times I_{m/(k-2)}\times\cdots\times  I_{m/(k-2)}$ is an  $m\times (m/(k-2))^{k-2}$
simple matrix  and $\lim_{m\to\infty} \ncols{A}/m^{k-2}= \frac{1}{(k-2)^{k-2}}$. 
Note that  $\lim_{m\to\infty} f(m,k)/m^{k-2}= \frac{2}{(k-2)!}$  leaving a wide gap for exact bounds for cases in \rf{6boundary}.

Coding theory helps to give a better construction for $\Av(m,F_{i,k})$ than the one given by the  product construction. One has to observe that $I_2\times \mathbf{1}_{k-2}\prec F_{i,k}$  for all $1\le i\le 6$. Let $A(m,4,k-1)$ be a constant weight code of minimum distance 4, length $m$ and weight $k-1$. Thus, if $A$ is an $m$-rowed matrix consisting of all columns with at most $k-2$ 1's and columns of codewords of $A(m,4,k-1)$, then $A\in\Av(m,I_2\times \mathbf{1}_{k-2})\subset\Av(m,F_{i,k})$. Indeed, $A$ does not have two columns with $k-1$ 1's on the same $k$-set of rows. By a theorem of Graham and Sloane \cite{GS80} this gives the lower bound
$\forb(m,F_{i,k})\ge\frac{1}{m}\binom{m}{k-1}+\sum_{i=0}^{k-2}\binom{m}{i}$. Thus, we have applying Theorem~\ref{simple,k-2}
\begin{E}\label{eq:code-lwbd}
    \frac{1}{k-1}\binom{m-1}{k-2}+\sum_{i=0}^{k-2}\binom{m}{i} \le \forb(m,F_{i,k})\le \binom{m-1}{k-2}+\sum_{i=0}^{k-2}\binom{m}{i}.
\end{E}
The lower bound of (\ref{eq:code-lwbd}) divided by $m^{k-2}$ decreases the wide gap mentioned above to 
\begin{E}\label{eq:codebd}
    \frac{1}{(k-1)!}+\frac{1}{(k-2)!}\le \lim_{m\to\infty}\frac{1}{m^{k-2}}\forb(m,F_{i,k})\le \frac{2}{(k-2)!}.
\end{E}
In fact, for $k=4,5$, Dukes \cite{Du15}  gives better lower bounds using nested block designs. Let $d(m,k)=\frac{1}{m^{k-2}}\forb(m,I_2\times \mathbf{1}_{k-2})$. His constructions give $d(m,4)\ge 0.6909$ and $d(m,5)\ge 0.25138$, while (\ref{eq:codebd}) results in only $d(m,4)\ge \frac{2}{3}$ and $d(m,5)\ge \frac{5}{24}\approx 0.2083$. For $k>5$, the difficult part of Duke's method is  finding good $t-(n, \ell, 1)$ packings where $n$ is not too much larger than $\ell$ for $t>3$. 
These results highlight the difficulty of finding the best, or indeed good constructions.

There are other functions that satisfy the recurrence \rf{recurrence} with different base cases including of course the binomial coefficients (with base cases $\binom{m}{k}=0$ for $m<k$ and $\binom{m}{0}=1$).  One standard proof for $\forb(m,K_k)$ uses the recurrence
  $\forb(m,K_k)=\forb(m-1,K_k)+\forb(m-1,K_{k-1})$ but in this case the base cases are $\forb(m,K_1)=1$ and $\forb(1,K_k)=2$ for $k\ge 2$.

The recurrence \rf{recurrence} appears in Geometry counting arguments \cite{GW} (see their Fact 2). For example $\forb(m,K_k)$ gives a bound on the number of regions in ${\bf R}^{k-1}$ when divided by $m$ hyperplanes.  The $k$ in $f(m,k)$ is typically the dimension while $m$ is the number of some geometric objects.  There are discussions of VC-dimension and Geometry not involving \rf{recurrence}.

\section{Constructions}\label{Balin}

  When we use induction to establish  the bound \trf{simple,k-2}, we forbid the three 2-rowed matrices
$${\cal F}=\left\{\left[\begin{array}{ccc}0&1&0\\ 0&0&1\\ \end{array}\right], \left[\begin{array}{ccc}0&1&1\\ 0&0&1\\ \end{array}\right],
\left[\begin{array}{ccc}0&1&1\\ 1&0&1\\ \end{array}\right]\right\},$$
for which $\forb(m,{\cal F})=2$.  This is the case $k=2$ in our recurrence.
As noted in \cite{survey} when $|{\cal F}|=1$, there are  6  possible column maximal $k$-rowed simple matrices satisfying the hypothesis of 
\trf{simple,k-2}.   
To obtain the list of \rf{6boundary},  we consider the pair of rows $i_1,i_2$ which has no 
$K_2^0=\left[\linelessfrac{0}{0}\right]$, the pair of rows  $i_3,i_4$ which has no $K_2^1=\left[\linelessfrac{1}{1}\right]$, and the pair of rows $i_5,i_6$ which has no $K_2^1=I_2$.
$F_{1,k}$ arises from the case  where the  three pairs of rows are distinct but whose union is just three rows. Note that $F_{1,3}=F_1$ in \cite{survey}.   The rows 
$\{i_1,i_2,i_3,i_4,i_5,i_6\}$ being distinct yields $F_{6,k}$. All cases are one of the six. For example if $i_1=i_3$ and $i_2=i_4$ while $i_5,i_6$
are disjoint from those, this yields $F_{2,k}$. 
If $i_1=i_3$ and $i_2=i_4$ while $i_2=i_5$ and $i_6$
is disjoint,  then the resulting columns are a subconfiguration of $F_{1,k}$.

 $$
 \hbox{Let }F_3=\left[\begin{array}{cc}1&1\\ 1&0\\ 0&1\\ \end{array}\right], 
F_4=\left[\begin{array}{cccc}1&1&1\\ 1&0&0\\ 0&1&0\\ \end{array}\right],
F_5=\left[\begin{array}{cccc}1&1&0\\ 1&0&1\\ 0&0&0\\ \end{array}\right]. 
$$
 Note that $F_3\prec F_4\prec F_{1,3}$ and $F_3^c\prec F_5\prec F_{1,3}$. In \cite{survey}, $F_{1,3}$ is called $F_1$.
 
 \begin{thm}\cite{survey}\label{F1} $\forb(m,F_{1,3})=2m=\forb(m,F_3)
 =\forb(m,F_3^c)$.  \end{thm}

\proof The upper bound $\forb(m,F_{1,3})\le 2m$ is \trf{simple,k-2}.  We can form an $m$-rowed simple $m\times 2m$ matrix  
$A=[\0\, |\,I_m\,|\,T_m^c\backslash\0_{m-1}\1_{1}]$ 
satisfying $F_3\not\prec A$. \qed

\vskip 10pt
There is more to be said in consideration of $F_{1,3}$'s role in \trf{F1}. 
We  explore the `what is missing' idea (\cite{survey}). One way to avoid $F_{1,3}$ is to have for each triple of rows $i,j,k$ (with $i<j<k$)

 $$\begin{array}{c@{}}\\ i\\ j\\ k\\ \end{array}\begin{array}{c}\hbox{no}\\ \left[\begin{array}{c}1\\ 0\\ 1\\  \end{array}\right]\end{array}\hbox{ and }
\begin{array}{c@{}}\\ i\\ j\\ k\\ \end{array} \begin{array}{c}\hbox{no}\\ \left[\begin{array}{c} 0\\ 1\\ 0\\  \end{array}\right]\end{array}$$ which yields an $m\times 2m$ simple matrix $ [T_m\,T_m^c]$ avoiding $F_{1,3}$. 
 This particular construction has $F_{1,3}\not\prec [T\,T^c]$ and it usefully generalizes.

\begin{lemma} Let $k\ge 3$. Let $A(k)$ be the $m$-rowed matrix  of all columns such that with  for each $k$-tuple $i_1<i_2<\cdots <i_k$, the columns satisfy
\begin{E}\begin{array}{c@{}}\\ i_1\\ i_2\\i_3\\ \vdots\\  i_k\\ \end{array}\begin{array}{c}\hbox{no}\\ \left[\begin{array}{c}1\\ 0\\ 1\\  \vdots\\ \\ \end{array}\right]\end{array}\hbox{ and }
\begin{array}{c@{}}\\ i_1\\ i_2\\i_3\\ \vdots\\  i_k\\ \end{array} \begin{array}{c}\hbox{no}\\ \left[\begin{array}{c} 0\\ 1\\ 0\\ \vdots\\ \\ \end{array}\right]\end{array}.\label{twomissing}\end{E} Then what is missing in each $k$-tuple of rows of  $A(k)$  is a pair of two complementary columns where for $k$ even, the complementary columns have $k/2$ 1's and for $k$ odd, the columns have  $\lceil k/2\rceil$ 1's and $\lfloor k/2\rfloor$ 1's.
Also $\ncols{A(k)}=f(m,k)$. \label{A(k)}\end{lemma}

\proof  The observation for what is missing follows from \rf{twomissing}.  Computing $\ncols{A(k)}$ is a little more work.
In a column, define a \emph{transition} at row $i$ in the adjacent rows $i,i+1$ if 
$$\begin{array}{c@{}}i\\ i+1\\ \end{array}\left[\begin{array}{c}0\\ 1\\ \end{array}\right]\hbox{ or }
\begin{array}{c@{}}i\\ i+1\\ \end{array}\left[\begin{array}{c}1\\ 0\\ \end{array}\right].$$
The `what is missing' conditions \rf{twomissing} force there to be at most $k-2$ transitions in any column.  If we choose $\binom{m-1}{k-2}$ rows 
$j_1<j_2<\cdots <j_{k-2}$ to be the rows for transitions then there are exactly  two columns $\alpha,\beta$ on $m$ rows with this pattern. Define
$\alpha$ as
 $\alpha|_{j_1}=0$, $\alpha|_{{j_1+1}}=1$ (a transition at row $j_1$) and 
 $\alpha|_{j_2}=1$, $\alpha|_{{j_2+1}}=0$ (a transition at row $j_2$ )and  
 $\alpha|_{j_3}=0$, $\alpha|_{{j_3+1}}=1$ (a transition at row $j_3$ )
 etc.  They must alternate else there will be additional transitions.  The remaining entries  of $\alpha$ are forced if there no other transitions.  Column $\beta$ has $\beta=\alpha^c$.  
 Thus the number of columns with $k-2$ transitions is $2\binom{m-1}{k-2}$.  Similarly the number of columns with (exactly) $t$ transitions
is $2\binom{m-1}{t}$ and so the number of columns satisfying \rf{twomissing} is
$$\sum_{i=0}^{k-2}2\binom{m-1}{i}=f(m,k),$$
establishing $\ncols{A(k)}=f(m,k)$ using \rf{f(m,k)}. \qed
\vskip 10pt
 As an example consider  $k=3$ which has $\binom{k}{2}=3$ complementary pairs of sum 1 and 2 .  We note that $F_3$, $F_4$ both have the property of having one  representative of each complementary pair and so $F_3,F_4\not\prec A(3)$.  Thus $\forb(m,F_3)=\forb(m,F_4)=2m$ by \trf{simple,k-2} for $k=3$ and the construction we have now provided.   Note that $F_3,F_4\prec F_{1,3}$.
 Recall that we do allow row and column permutations after $A(k)$ has been formed.
\vskip 10pt
\noindent {\bf Proof of \trf{mainthm}:}
For $k$ even, there are $\frac{1}{2}\binom{k}{k/2}$ complementary pairs of columns of sum $k/2$.  For $k$ odd there are 
$\binom{k}{\lceil k/2\rceil}$ complementary pairs of a column of sum $\lceil k/2\rceil$ and a column of sum $\lfloor k/2\rfloor$. 
    If a $k$-rowed $F$ has one representative of each such complementary pair, then $F\not\prec A(k)$. We use this to consider $A(k)$

For $k=4$, there are $\frac{1}{2}\binom{k}{2}=3$ complementary pairs.  
We consider the 4-rowed  matrices $F_{1,4}$ and $F_{3,4}$. 

$$\hbox{Let } 
F_5=\left[\begin{array}{ccc}1&1&0\\ 1&0&1\\ 0&1&1\\ 0&0&0\\ \end{array}\right],
F_6=\left[\begin{array}{ccc}1&0&0\\ 0&1&0\\ 0&0&1\\ 1&1&1\\ \end{array}\right] $$
 Both $F_5,F_6\prec F_{1,4},F_{3,4}$ and each have one  representative of each complementary pair of sum 2 and so $F_5,F_6\not\prec A(4)$. 
 Thus $\forb(m,F_5),\forb(m,F_6)\ge f(m,4)$.
Given $F_5,F_6\prec F_{1,4},F_{3,4}$, we deduce $\forb(m,F_5)=\forb(m,F_6)=f(m,4)=2\binom{m}{2}+2\binom{m}{0}$ by \trf{simple,k-2}.  
  Also for any $F$ with $F_5\prec F\prec F_{1,4}$ or   $F_6\prec F\prec F_{1,4}$ or $F_5\prec F\prec F_{3,4}$ or   $F_6\prec F\prec F_{3,4}$, we have 
  $\forb(m,F)=f(m,4)$.

Consider the  matrices $F_{3,k}$ for $k=5,6$. On $k$ rows, these matrices have (at least) one representative of each  pair of complementary columns  of column sum 2,3 for $k=5$ and column sum $3$ for $k=6$.  


For $k=5$, for each of the 10 columns of sum 2, there is a complementary column of sum 3.
We need to verify that $F_{3,5}$ has a representative of each such complementary pair.
In particular $F_{3,5}$ has 6 columns of sum 2 and 6 columns of sum 3 with a representative of each such complementary pair
(there are exactly two  complementary pairs with both present among those 12 columns).  The matrix $A(5)$ from \lrf{A(k)}  avoids such a complementary pair on each $5$-set of rows and $F_{3,5}\not\prec A(5)$.  
Thus $\forb(m,F_{3,5})=f(m,5)$.  But we can consider the minimal ways to have a representative of each such complementary pairs (of columns of sum 2,3) while remaining in $F_{3,5}$, yielding 4 $5\times 10$ matrices $F$ for which each has $\forb(m,F)=f(m,5)$.

$$\left[\begin{array}{cccccccccc}
0& 0& 0& 0& 0& 0& 0& 1& 1& 1\\
0& 0& 0& 0& 1& 1& 1& 0& 0& 0\\
0& 1& 1& 1& 0& 0& 1& 1& 0& 0\\
1& 0& 1& 1& 1& 1& 1& 1& 1& 1\\
1& 1& 0& 1& 0& 1& 0& 0& 0& 1\\
 \end{array}\right],\quad
\left[\begin{array}{cccccccccc}
0& 0& 0& 0& 0& 0& 0& 1& 1& 0\\
0& 0& 0& 0& 1& 1& 1& 0& 0& 1\\
0& 1& 1& 1& 0& 0& 1& 1& 0& 1\\
1& 0& 1& 1& 1& 1& 1& 1& 1& 0\\
1& 1& 0& 1& 0& 1& 0& 0& 0& 0\\
 \end{array}\right],$$
$$\left[\begin{array}{cccccccccc}
0& 0& 0& 0& 0& 0& 0& 1& 0& 0\\
0& 0& 0& 0& 1& 1& 1& 0& 1& 1\\
0& 1& 1& 1& 0& 0& 1& 1& 1& 1\\
1& 0& 1& 1& 1& 1& 1& 1& 0& 0\\
1& 1& 0& 1& 0& 1& 0& 0& 0& 1\\
 \end{array}\right],\quad
\left[\begin{array}{cccccccccc}
0& 0& 0& 0& 0& 0& 0& 1& 0& 1\\
0& 0& 0& 0& 1& 1& 1& 0& 1& 0\\
0& 1& 1& 1& 0& 0& 1& 1& 1& 0\\
1& 0& 1& 1& 1& 1& 1& 1& 0& 1\\
1& 1& 0& 1& 0& 1& 0& 0& 1& 1\\
  \end{array}\right]$$

For $k=6$, there are $\frac{1}{2}\binom{6}{3}=10$ complementary pairs of columns of sum 3.  In $F_{3,6}$, there are exactly 12 columns of sum 3 with a representative  of each of the 10 complementary pairs of sum 3.  If a matrix avoids a complementary pair of columns of sum 3 on each 6-set of rows, then it avoids $F_{3,6}$. Thus $F_{3,6}\not\prec A(6)$ so  $\forb(m,F_{3,6})=f(m,6)$.  But we can consider the minimal ways to have a representative of each of the 10 complementary pairs of columns of sum 3 while remaining in $F_{3,6}$, yielding 4 $6\times 10$ matrices $F$ for which each has $\forb(m,F)=f(m,6)$.

$$\left[\begin{array}{cccccccccc}
0& 0& 0& 0& 0& 0& 0& 1& 1& 1\\
0& 0& 0& 0& 1& 1&  1& 0& 0& 0\\
0& 1& 1&  1& 0& 0& 1& 1& 0& 0\\
1& 0& 1& 1&  1& 1& 1&  1& 1& 1\\
1& 1&  1& 0& 1& 0& 0& 0& 1& 0\\
1& 1& 0& 1& 0& 1& 0& 0& 0& 1\\
 \end{array}\right],\quad
\left[\begin{array}{cccccccccc}
0& 0& 0& 0& 0& 0& 0& 1& 1& 0\\
0& 0& 0& 0& 1& 1&  1& 0& 0& 0\\
0& 1& 1&  1& 0& 0& 1& 1& 0& 1\\
1& 0& 1& 1&  1& 1& 1&  1& 1& 0\\
1& 1&  1& 0& 1& 0& 0& 0& 1& 1\\
1& 1& 0& 1& 0& 1& 0& 0& 0& 0\\
  \end{array}\right],$$

$$\left[\begin{array}{cccccccccc}
0& 0& 0& 0& 0& 0& 0& 1& 1& 0\\
0& 0& 0& 0& 1& 1&  1& 0& 0& 1\\
0& 1& 1&  1& 0& 0& 1& 1& 0& 1\\
1& 0& 1& 1&  1& 1& 1&  1& 1& 0\\
1& 1&  1& 0& 1& 0& 0& 0& 0& 0\\
1& 1& 0& 1& 0& 1& 0& 0& 1& 1\\
\end{array}\right],
\quad
\left[\begin{array}{cccccccccc}
0& 0& 0& 0& 0& 0& 0& 1& 0& 0\\
0& 0& 0& 0& 1& 1&  1& 0& 1& 1\\
0& 1& 1&  1& 0& 0& 1& 1& 1& 1\\
1& 0& 1& 1&  1& 1& 1&  1& 0& 0\\
1& 1&  1& 0& 1& 0& 0& 0& 1& 0\\
1& 1& 0& 1& 0& 1& 0& 0& 0& 1\\
  \end{array}\right]$$
\qed\vskip 10pt

For $k=8$, there are representatives of only 34 complementary pairs of columns of sum 4 in $F_{3,8}$ when there are 35 complementary pairs in total.  Thus  $A(8)$ may not avoid $F_{3,8}$ although ${A(8)}$ is a simple matrix achieving the bound $\ncols{A(8)}= f(m,8)$.

\end{document}